\documentclass{amsart}

\usepackage{amssymb,bm,graphicx}
\usepackage{hyperref}

\newtheorem{thm}{Theorem}
\newtheorem{lem}[thm]{Lemma}

\newcommand{\gauss}[2]{\genfrac{[}{]}{0pt}{}{#1}{#2}}
\renewcommand{\qed}{\hfill\resizebox{1.2ex}{.7em}{$\blacksquare$}}

\begin{document}

\title{The Erd\H{o}s--Ko--Rado theorem for twisted Grassmann graphs}
\author{Hajime Tanaka}
\address{Department of Mathematics, University of Wisconsin, 480 Lincoln Drive, Madison, WI 53706, U.S.A.}
\curraddr{Graduate School of Information Sciences, Tohoku University, 6-3-09 Aramaki-Aza-Aoba, Aoba-ku, Sendai 980-8579, Japan}
\email{htanaka@math.is.tohoku.ac.jp}
 % 05E30: Algebraic combinatorics; 05D05: Extremal set theory
\begin{abstract}
We present a ``modern'' approach to the Erd\H{o}s--Ko--Rado theorem for $Q$-polynomial distance-regular graphs and apply it to the twisted Grassmann graphs discovered in 2005 by van Dam and Koolen.
\end{abstract}

\maketitle

%%%%%%%%%%%%%%%%%%%%%%%%%%%%%%%%%%%%%%%%%%%
%%%%%%%%%%%%%%%%%%%%%%%%%%%%%%%%%%%%%%%%%%%
\section{Introduction}
The 1961 theorem of Erd\H{o}s, Ko and Rado \cite{EKR1961QJMO} asserts that the largest possible families $Y$ of $d$-subsets of a $v$-set such that $|x\cap y|\geqslant t$ for all $x,y\in Y$ where $v>(t+1)(d-t+1)$ are the families of all $d$-subsets containing some fixed $t$-subset.
In fact, the exact bound $v>(t+1)(d-t+1)$ was obtained later by Wilson \cite{Wilson1984C} as an application of Delsarte's linear programming method \cite{Delsarte1973PRRS}.
It is natural to think of this theorem as a result about (vertex) subsets of the Johnson graphs $J(v,d)$, and analogous theorems are known for several other families of distance-regular graphs, e.g., Hamming graphs $H(d,q)$ $(q\geqslant t+2)$ \cite{Moon1982JCTA}, Grassmann graphs $J_q(v,d)$ $(v\geqslant 2d)$ \cite{Hsieh1975DM,FW1986JCTA,Fu1999DM,Tanaka2006JCTA}, bilinear forms graphs $\mathrm{Bil}_q(d,e)$ $(d\leqslant e)$ \cite{Huang1987DM,Fu1999DM,Tanaka2006JCTA}.

In this note, we first distill common algebraic techniques found in some of the proofs of these ``Erd\H{o}s--Ko--Rado theorems'' into a unified approach for general $Q$-polynomial distance-regular graphs $\Gamma$.\footnote{$Q$-polynomial distance-regular graphs are thought of as finite/combinatorial analogues of symmetric spaces of rank one; see \cite[pp.~311--312]{BI1984B}.}
Our approach is also ``modern'' in the sense that it is based on and motivated by the theory of two parameters, \emph{width} $w$ and \emph{dual width} $w^*$, of a subset $Y$ of $\Gamma$ introduced in 2003 by Brouwer et al.~\cite{BGKM2003JCTA}.
In this setting, the ``$t$-intersecting'' condition amounts to requiring $w\leqslant d-t$ where $d$ is the diameter of $\Gamma$, and we shall view the Erd\H{o}s--Ko--Rado theorem as characterizing those subsets $Y$ with $w=d-t$ and $w^*=t$ by their sizes among all $t$-intersecting families.
There are two steps involved: (1) construction of a specific feasible solution to the dual of a linear programming problem; (2) classification of the \emph{descendents} \cite{Tanaka2010pre} of $\Gamma$, i.e., those subsets having the property $w+w^*=d$.
We demonstrate this approach by deriving the Erd\H{o}s--Ko--Rado theorem for the \emph{twisted Grassmann graphs} $\tilde{J}_q(2d+1,d)$ discovered in 2005 by van~Dam and Koolen \cite{DK2005IM}.

%%%%%%%%%%%%%%%%%%%%%%%%%%%%%%%%%%%%%%%%%%%
%%%%%%%%%%%%%%%%%%%%%%%%%%%%%%%%%%%%%%%%%%%
\section{A ``modern'' approach to the Erd\H{o}s--Ko--Rado theorem for $Q$-polynomial distance-regular graphs}\label{sec: modern approach}

Let $\Gamma=(X,R)$ be a finite connected simple graph with diameter $d$ and path-length distance $\partial$, and $\mathbb{R}^{X\times X}$ the set of real matrices with rows and columns indexed by $X$.
For each $i$ $(0\leqslant i\leqslant d)$, let $A_i\in\mathbb{R}^{X\times X}$ be the adjacency matrix of the distance-$i$ graph $\Gamma_i$ of $\Gamma$, so $A_0=I$ and $\sum_{i=0}^dA_i=J$, the all ones matrix.
We say $\Gamma$ is \emph{distance-regular} if $\bm{A}:=\mathrm{span}\{A_0,A_1,\dots,A_d\}$ is closed under ordinary matrix multiplication; or equivalently, $\bm{A}$ is a (commutative) algebra.
(The reader is referred to \cite{BI1984B,BCN1989B,Godsil1993B} for background material on distance-regular graphs.)
Throughout this note, suppose $\Gamma$ is distance-regular.
We call $\bm{A}$ the \emph{Bose--Mesner algebra} of $\Gamma$.
It is semisimple (as it is closed under transposition) and therefore has a basis $\{E_i\}_{i=0}^d$ consisting of the  primitive idempotents;
we always set $E_0=|X|^{-1}J$.
Note that $\bm{A}$ is also closed under entrywise multiplication, denoted $\circ$.
We shall assume $\Gamma$ is $Q$-\emph{polynomial} with respect to the ordering $\{E_i\}_{i=0}^d$, i.e., $E_1\circ E_i$ is a linear combination of $E_{i-1},E_i,E_{i+1}$ with nonzero coefficients for $E_{i-1},E_{i+1}$ $(0\leqslant i\leqslant d)$, where $E_{-1}=E_{d+1}=0$.
Let $Q=(Q_{ij})_{0\leqslant i,j\leqslant d}$ be the \emph{second eigenmatrix} of $\Gamma$:
\begin{equation*}
	E_j=\frac{1}{|X|}\sum_{i=0}^dQ_{ij}A_i \quad (0\leqslant j\leqslant d).
\end{equation*}

Let $Y$ be a nonempty subset of $X$ and $\chi\in\mathbb{R}^X$ its (column) characteristic vector.
Brouwer et al.~\cite{BGKM2003JCTA} defined the \emph{width} $w$ and \emph{dual width} $w^*$ of $Y$ as follows:
\begin{equation*}
	w=\max\{i:\chi^{\mathsf{T}}A_i\chi\ne 0\}, \quad w^*=\max\{i:\chi^{\mathsf{T}}E_i\chi\ne 0\}.
\end{equation*}
They showed (among other results) that
\begin{equation}\label{fundamental inequality}
	w+w^*\geqslant d.
\end{equation}
We call $Y$ a \emph{descendent} \cite{Tanaka2010pre} of $\Gamma$ if $w+w^*=d$.
It should be remarked that every descendent is a so-called completely regular code (cf.~\cite{KLM2010P}), and that the induced subgraph is a $Q$-polynomial distance-regular graph provided it is connected; see \cite[Theorems 1--3]{BGKM2003JCTA}.
See also \cite{Tanaka2010pre} for more information on descendents.

Now fix an integer $t$ $(0<t<d)$ and suppose $w\leqslant d-t$; in other words, $Y$ is ``$t$-intersecting''.
We recall the inner distribution $\bm{e}=(e_0,e_1,\dots,e_d)$ of $Y$:
\begin{equation*}
	e_i=\frac{1}{|Y|}\chi^{\mathsf{T}}A_i\chi, \quad (\bm{e}Q)_i=\frac{|X|}{|Y|}\chi^{\mathsf{T}}E_i\chi \quad (0\leqslant i\leqslant d).
\end{equation*}
It follows that $|Y|=(\bm{e}Q)_0$ and
\begin{gather*}
	e_0=1, \quad e_1\geqslant 0,\dots,e_{d-t}\geqslant 0,\quad e_{d-t+1}=\dots=e_d=0, \\
	(\bm{e}Q)_1\geqslant 0,\dots,(\bm{e}Q)_d\geqslant 0.
\end{gather*}
(Observe that the $E_i$ are positive semidefinite.)
Following \cite{Delsarte1973PRRS}, we view these as a linear programming maximization problem.
A vector $\bm{f}=(f_0,f_1,\dots,f_d)$ satisfying \eqref{constraint 1}, \eqref{constraint 2} below gives a feasible solution to its dual problem:
\begin{gather}
	f_0=1, \quad f_1=\dots=f_t=0, \quad f_{t+1}>0,\dots,f_d>0, \label{constraint 1} \\
	(\bm{f}Q^{\mathsf{T}})_1=\dots=(\bm{f}Q^{\mathsf{T}})_{d-t}=0. \label{constraint 2}
\end{gather}
Indeed, we have
\begin{equation*}
	|Y|=(\bm{e}Q)_0\leqslant\bm{e}Q\bm{f}^{\mathsf{T}}=(\bm{f}Q^{\mathsf{T}})_0
\end{equation*}
with equality if and only if $(\bm{e}Q)_{t+1}=\dots=(\bm{e}Q)_d=0$, i.e., $w^*\leqslant t$.
By virtue of \eqref{fundamental inequality}, it follows that

\begin{lem}\label{EKR for general Gamma}
Let $Y$ be a nonempty subset of $X$ with $w\leqslant d-t$.
Suppose there is a vector $\bm{f}=(f_0,f_1,\dots,f_d)$ satisfying \eqref{constraint 1}, \eqref{constraint 2}.
Then $|Y|\leqslant(\bm{f}Q^{\mathsf{T}})_0$, and equality holds if and only if $Y$ is a descendent of $\Gamma$ with $w=d-t$ and $w^*=t$. \qed
\end{lem}

The vector $\bm{f}$ above is of independent interest from the point of view of \emph{Leonard systems}\footnote{Leonard systems provide a linear algebraic framework characterizing the terminating branch of the Askey scheme \cite{KLS2010B} of (basic) hypergeometric orthogonal polynomials.} \cite{Terwilliger2006N} and will be discussed in detail in a future paper.
Here we mention that $\bm{f}$ can be found for the following graphs:
\smallskip
\begin{center}
\begin{tabular}{lc}
\hline\hline
\multicolumn{1}{c}{$\Gamma$} & $(\bm{f}Q^{\mathsf{T}})_0$ \\
\hline
$J(v,d)$ ($v>(t+1)(d-t+1)$) & $\binom{v-t}{d-t}$ \\
%$H(d,q)$ ($q\geqslant d$; or $q=d-1$, $t\ne d-2$) & $q^{d-t}$ \\
$H(d,q)$ ($t=d-1$; or $q\geqslant d$; or $q=d-1$, $t<d-2$) & $q^{d-t}$ \\
$J_q(v,d)$ ($v\geqslant 2d$) & $\gauss{v-t}{d-t}_q$ \\
$\mathrm{Bil}_q(d,e)$ ($d\leqslant e$) & $q^{(d-t)e}$ \\
\hline\hline
\end{tabular}
\end{center}
\smallskip
For $\Gamma=J(v,d)$ or $J_q(v,d)$ (with $v,d$ as in the table), Wilson and Frankl \cite{Wilson1984C,FW1986JCTA} constructed a matrix $B\in\bm{A}$ such that (i) $B_{xy}=0$ if $\partial(x,y)\leqslant d-t$; (ii) $B+I-\gauss{v-t}{d-t}^{-1}J$ is positive semidefinite and its $i^{\text{th}}$ eigenvalue $\lambda_i$ is positive precisely when $t+1\leqslant i\leqslant d$, where we interpret $\gauss{m}{n}$ as $\binom{m}{n}$ for $J(v,d)$ and $\gauss{m}{n}_q$ for $J_q(v,d)$.
We define $\bm{f}$ by $f_0=1$, $f_1=\dots=f_t=0$, and $f_i=\gauss{v-t}{d-t}\gauss{v}{d}^{-1}\lambda_i$ for $t+1\leqslant i\leqslant d$.
For $\Gamma=\mathrm{Bil}_q(d,e)$ $(d\leqslant e)$, Delsarte \cite{Delsarte1978JCTA} constructed a \emph{Singleton system}, i.e., a subset whose inner distribution $\bm{e}'=(e_0',e_1',\dots,e_d')$ satisfies $e_1'=\dots=e_t'=0$ and $(\bm{e}'Q)_1=\dots=(\bm{e}'Q)_{d-t}=0$.
It follows that $e_{t+1}',\dots,e_d'$ are positive; see \cite[\S 4]{Tanaka2006JCTA}.
We define $\bm{f}=\bm{e}'\cdot\mathrm{diag}(k_0,k_1,\dots,k_d)^{-1}$ where $k_i$ is the valency of $\Gamma_i$ $(0\leqslant i\leqslant d)$.
For $\Gamma=H(d,q)$, a subset having the above properties is known as an MDS code \cite[Chapter 11]{MS1977B}.
MDS codes may not exist for some $d,q,t$, but still $\bm{e}'$ makes sense and is uniquely determined.
If $t=d-1$ or $q\geqslant d$, or if $q=d-1$ and $t<d-2$, then it follows that $e_{t+1}',\dots,e_d'$ are positive; see e.g., \cite[Appendix]{EGS2011IEEE}.
We again define $\bm{f}=\bm{e}'\cdot\mathrm{diag}(k_0,k_1,\dots,k_d)^{-1}$.

Concerning the conclusion of Lemma \ref{EKR for general Gamma}, the classification of descendents has been done for the 15 known infinite families of $Q$-polynomial distance-regular graphs with so-called classical parameters and with unbounded diameter, including the above 4 families; see \cite{BGKM2003JCTA,Tanaka2006JCTA,Tanaka2010pre}.
Moon \cite{Moon1982JCTA} showed that the upper bound $q^{d-t}$ for $H(d,q)$ and the characterization of its descendents as optimal intersecting families are valid under the (in general) weaker assumption $q\geqslant t+2$.
Dual polar graphs discussed in \cite{Tanaka2006JCTA} do not always possess $\bm{f}$ even for the case $t=1$ \cite{Stanton1980SIAMb}; see \cite{PSV2011JCTA}, however, for a description of optimal $1$-intersecting families.

%%%%%%%%%%%%%%%%%%%%%%%%%%%%%%%%%%%%%%%%%%%
%%%%%%%%%%%%%%%%%%%%%%%%%%%%%%%%%%%%%%%%%%%
\section{The Erd\H{o}s--Ko--Rado theorem for twisted Grassmann graphs}

Let $q$ be a prime power and fix a hyperplane $H$ of $\mathbb{F}_q^{2d+1}$.
Let $X_1$ be the set of $(d+1)$-dimensional subspaces of $\mathbb{F}_q^{2d+1}$ not contained in $H$, and $X_2$ the set of $(d-1)$-dimensional subspaces of $H$.
The \emph{twisted Grassmann graph} $\Gamma=\tilde{J}_q(2d+1,d)$ \cite{DK2005IM} has vertex set $X=X_1\cup X_2$, and two vertices $x,y\in X$ are adjacent if $\dim x+\dim y-2\dim x\cap y=2$.
It has the same parameters (i.e., the structure constants of $\bm{A}$) as $J_q(2d+1,d)$.
The twisted Grassmann graphs provide the first known family of non-vertex-transitive distance-regular graphs with unbounded diameter.
See \cite{FKT2006IIG,BFK2009EJC,MT2009pre} for more information.

The Erd\H{o}s--Ko--Rado theorem for $\tilde{J}_q(2d+1,d)$ can now be rapidly obtained.
Note that $J_q(2d+1,d)$ and $\tilde{J}_q(2d+1,d)$ share the same $Q$.
Hence we may use the vector $\bm{f}$ for $J_q(2d+1,d)$ constructed in \S \ref{sec: modern approach}, and Lemma \ref{EKR for general Gamma} applies.
The descendents of $\tilde{J}_q(2d+1,d)$ have recently been classified by the author \cite[Theorem 8.20]{Tanaka2010pre}.
To summarize:

\begin{thm}
Let $Y$ be a nonempty subset of $\tilde{J}_q(2d+1,d)$ with width $w\leqslant d-t$, where $0<t<d$.
	Then $|Y|\leqslant\gauss{2d+1-t}{d-t}_q$, and equality holds if and only if $Y=\{x\in X_2:u\subseteq x\}$ for some subspace $u$ of $H$ with $\dim u=t-1$. \qed
\end{thm}

%%%%%%%%%%%%%%%%%%%%%%%%%%%%%%%%%%%%%%%%%%%
%%%%%%%%%%%%%%%%%%%%%%%%%%%%%%%%%%%%%%%%%%%
\section*{Acknowledgements}

The author would like to thank the Department of Mathematics at the University of Wisconsin--Madison for its hospitality throughout the period in which this work was done.
Support from the JSPS Excellent Young Researchers Overseas Visit Program is also gratefully acknowledged.

\end{document}